\documentclass[preprint]{elsarticle}
\usepackage[margin=3cm]{geometry}
\usepackage{graphicx}
\usepackage{floatrow}
\usepackage{amsmath}

\begin{document}
\begin{frontmatter}
\title{\textit{dune-composites} -- A New Framework for
  High-Performance Finite Element
  Modelling of Laminates}

\author[bm]{Anne Reinarz}
\address[bm]{Department of Mathematical Sciences, University of Bath}
\author[e]{Tim Dodwell}
\address[e]{College of Engineering, Mathematics and Physical Sciences, University of Exeter}
\author[be]{Tim Fletcher}
\address[be]{Department of Mechanical Engineering, University of Bath}
\author[h]{Linus Seelinger}
\address[h]{Institute for Scientific Computing, University of Heidelberg}
\author[be]{Richard Butler}
\author[bm]{Robert Scheichl}

\begin{abstract}
Finite element (FE) analysis has the potential to offset much of the expensive
experimental testing currently required to certify aerospace laminates. However,
large numbers of degrees of freedom are necessary to model entire aircraft components
whilst accurately resolving micro-scale defects. The new module \textit{dune-composites},
implemented within DUNE by the authors, provides a tool to efficiently solve large-scale problems using novel
iterative solvers. The key innovation is a preconditioner that guarantees a constant number of iterations 
regardless of the problem size. Its robustness has been shown rigorously in Spillane et al. ({\em Numer.~Math.} {\bf 126}, 2014) for isotropic problems. For anisotropic problems in composites it is verified numerically
for the first time in this paper. The parallel implementation in DUNE scales almost optimally over thousands of cores.
To demonstrate this, we present an original numerical study, varying the shape of a localised wrinkle
and the effect this has on the strength of a curved laminate. This requires a high-fidelity mesh containing 
at least four layers of quadratic elements across each ply and interface layer, underlining the need for 
\textit{dune-composites}, which can achieve run times of just over 2 minutes on 2048 cores for realistic 
composites problems with 173 million degrees of freedom. 

\end{abstract}
 \end{frontmatter}
 
\section{Introduction}
Traditionally, aerospace structures have been largely manufactured by subtractive processes, such as machining,
which is typically applied to metals. Composite materials, which constitute over 50\% of recent aircraft construction, 
are created at the same time as the structure itself; typically by an additive layup process such as
automated fibre placement (AFP). This co-assembly of material and structure 
is achieved by sequential deposition of a number of fibrous layers which are either pre-impregnated or post-infused with resin, involving
hundreds of additive operations during the production of a large structural part.
Most of these processes take place before the resin is cured, when the influence
of temperature and pressure on the deposition, forming and curing of parts with 
complex geometry is not well understood. For this reason production is vulnerable
to defects arising from small process variations and achieving high-rate automation
of the process is limited.

The ability to manufacture composite material-structural 
systems with repeatable properties at a reasonable rate and competitive cost,
with the quality required for certification, is a major challenge in the global 
aerospace industry. A number of defect types can arise, and use of non-destructive
evaluation is limited in detecting these, hence conservative `knock-down' factors 
are applied to strength limits. These factors are currently obtained by 
expensive testing at all scales in the Test Pyramid, with large numbers of small 
coupon tests and fewer larger-scale structural tests. 
Large-scale tests, are extremely expensive and take 
place at a stage when it is difficult to make changes and improve designs, while
small-scale coupons do not represent the performance of the material at the 
structural scale where, for example, defects can be introduced as a result of 
complex forming processes.
Furthermore,
the boundary conditions and loading in a simple element are very different from the performance
of the material within the structure.

One loading of particular interest is  corner unfolding,
in which though thickness tensile loading acts to separate the layers as the corner radius
is increased under bending loads. These stresses act in the weak resin-dominated direction of the
laminate and can therefore lead to a limiting design case. The presence of defects such as mis-aligned
fibres in the form of out-of-arc wrinkles can have a significant importance. Such a reduction
in strength is referred to in manufacturing as a knock-down factor.

Mukhopadhyay et al have modelled the effect of defects on compressive \cite{muk1}
and tensile \cite{muk2} strength in flat laminates using 8-node, solid elements
and zero-thickness, 8-node, cohesive elements between plies. Damage modelling accounted 
for nonlinear shear in plies, 
transverse matrix cracking,
mixed mode delamination, tensile fibre fracture and fibre kinking. The 
minimum size of the FE mesh was one ply thickness (0.25 mm) in the vicinity of wrinkles and towards laminate edges.
The compressive strength predictions were in agreement to within 10\%
of experimental results and were able to pick up a mode switch from
fibre failure to delamination when defect  misalignment was above
$\sim 9^\circ$. In the tensile case, it was noted that wrinkles
act as local, through-thickness shear stress concentrators. 
For multidirectional laminates, the influence 
of the defect was exacerbated by edge effects.

We have previously shown that the application of a $3$ mm tough resin layer to the longitudinal
edges of corner unfolding specimens prevents the zero stress requirement along the edge, thus reducing the numerical
singularity at the edge. The finite element solution becomes more realistic for a given
mesh fidelity and failure prediction
methods become more accurate. Such an edge treatment can increase the strength by up to $20\%$ 
\cite{fletcher2016resin} and provides a more accurate representation of a long structural part.

There is a wide range of defects that can form in composite laminates and a 
detailed taxonomy of these is presented in \cite{potter}. 
In this paper we evaluate wrinkle defects, which are more likely to form 
in curved laminates as a result of consolidation onto a male tool \cite{dodwell}. 
They are also likely to form in curved laminates with tapered sections, which 
causes double curvature and makes AFP deposition challenging.  These wrinkles are 
important to the performance of curved laminates.

In this paper we present a new FE analysis tool \textit{dune-composites} for efficient, high-fidelity
modeling of laminated composite parts. We show, using a simple test case that the results of \textit{dune-composites}
match up well with the commercial software package ABAQUS in all six
stresses. For large scale problems, \textit{dune-composites} crucially
relies on robust iterative solvers for the resulting FE systems. 
To this end, we introduce the preconditioner GenEO \cite{geneo1,geneo2,geneoparallel}, which we have implemented within DUNE.
This preconditioner has previously been mathematically proven to be
robust for isotropic FE systems.
We show that these results extend to anisotropic problems and that the solver is suitable
for solving large composites problems. Further, we demonstrate its
parallel efficiency and its ability to scale to
thousands of compute cores, allowing the solution of the large problems with defects mentioned above.
We test this module by modeling the unfolding of a curved laminate
part containing manufacturing defects, for which a micron scale mesh is needed to accurately compute stresses. We show that 
\textit{dune-composites} is able to accurately predict damage initiation and that it does so
at a fraction of the computational cost required by ABAQUS. 


\section{Modelling approach}

In this section, we introduce the new high performance finite element module \textit{dune-composites},
and demonstrate its capability of efficiently and robustly tackling
large-scale simulations of composite structures. In the example simulations presented, the composite strength of pristine
and defected corner radii are accurately predicted. The analysis assumes standard anisotropic
3D linear elasticity and the failure is assessed using a quadratic damage onset
criterion for the initiation of delamination in \cite{Camanho}.
This allows the results to be benchmarked against existing numerical results and experiments, given 
in Fletcher et al. \cite{fletcher2016resin}. For this case, the numerical results
of this linear model show good agreement with experimental data.

The failure of curved laminates subjected to corner unfolding has been shown 
to be highly unstable with instantaneous propagation following initiation \cite{fletcher2016resin}.
Therefore we focus on capturing the initiation of failure accurately by using a high-fidelity
3D mesh, with $4$ or more elements through the thickness of each ply and interface layer.
Such an approach has been shown to give accurate results for simple flat coupons \cite{esquej}
and 2D models of a few plies \cite{2dply}. However, this level of fidelity is typically dropped when modelling 
larger, more complex parts, such as 3D curved laminates, due to computational limitations. 
Here, we do not use any of these more complex models in order to reduce the computational requirements.

More complex modelling approaches have been proposed for modelling of defects in the literature,
these include the use of composite shell elements \cite{muk2, muk1},
interface/cohesive elements \cite{hallettcohesive, hallett2008}
and higher-order continuum models \cite{fleck2004}. In particular, cohesive elements are commonly
used to capture propagation \cite{hallettcohesive},
but propagation is considered less important than initiation for the problem.
The formulation of such 
elements, whilst more complex, is possible within \textit{dune-composites}. The solution
strategy would require Newton iterations and path-following methods (available in the DUNE library).
However, even for these non-linear models of failure, computational cost 
is still dominated by the speed of solving a linearised system of equations for composite materials.
In this paper our focus is therefore on developing and implementing an efficient solution strategy 
for these linearised equations arising from such massive composite simulation.


\subsection{{\em dune-composites}: High performance FE modelling of large-scale composite structures}
DUNE (Distributed and Unified Numerics Environment) is an open source modular toolbox 
for solving partial differential equations (PDEs) with grid-based methods, 
such as the finite element method (FEM) \cite{dune2a,dune2b,dune1}. Written using modern C++
programming techniques,
the core modules of DUNE have been developed by mathematicians and computer scientists
to allow users to implement and use state-of-the-art mathematical methods across large
high performance parallel computing architectures. It is a generic
package that provides a user with the key ingredients for solving any FEM problem,
e.g. grid generation, different types of finite elements, quadrature rules
and a choice of off-the shelf solvers. 

Within this platform, we have developed a new module, \textit{dune-composites}, which
solves the linear elasticity equations with general, anisotropic
stiffness tensor, applicable for modelling
 composite structures.  This module provides an interfaces to handle composite
 applications which includes stacking sequences, complex part geometries and
 complex boundary conditions such as multi-point constraints, or periodic boundary conditions.
Further, we have implemented a new $20$-node 3D serendipity element
(with full integration) within \textit{dune-pdelab},
 which is not prone to shear locking and
 allows comparison with ABAQUS's C3D20R element. This element has degrees of freedom
 at the $8$ nodes of the element as well as on each of the $12$ edges.
 
 The main advantage of implementing our simulation tool within DUNE
 is that it allows us to exploit  developments in state of the art solvers
 such as Algebraic Multigrid (AMG), see \cite{dune3}, or to implement
 new ones. In particular, as part of our new developments, we
 implemented a novel, robust preconditioner, called GenEO
 \cite{geneo2,geneo1}, within DUNE.
Finally, we will show in Section \ref{sec-time} that DUNE allows for  highly parallelised efficiency
on hundreds of computer cores. This allows for the modelling of meshes fine enough to
resolve defects or the modelling of wide parts with sufficient accuracy.

\subsection{Verification of DUNE}
\label{sec-verify}


In this section we show that DUNE  produces stress results that are comparable to 
those produced by standard FE libraries such as ABAQUS. 
For our comparison we examine a corner unfolding test in DUNE 
and ABAQUS and compare the solutions as well as the run times. We
initially use a small computationally inexpensive problem for
comparisons with ABAQUS and move to a larger problem to demonstrate the 
parallel efficiency of DUNE.

\begin{figure}[t]
 \centering
 \includegraphics[scale=0.9]{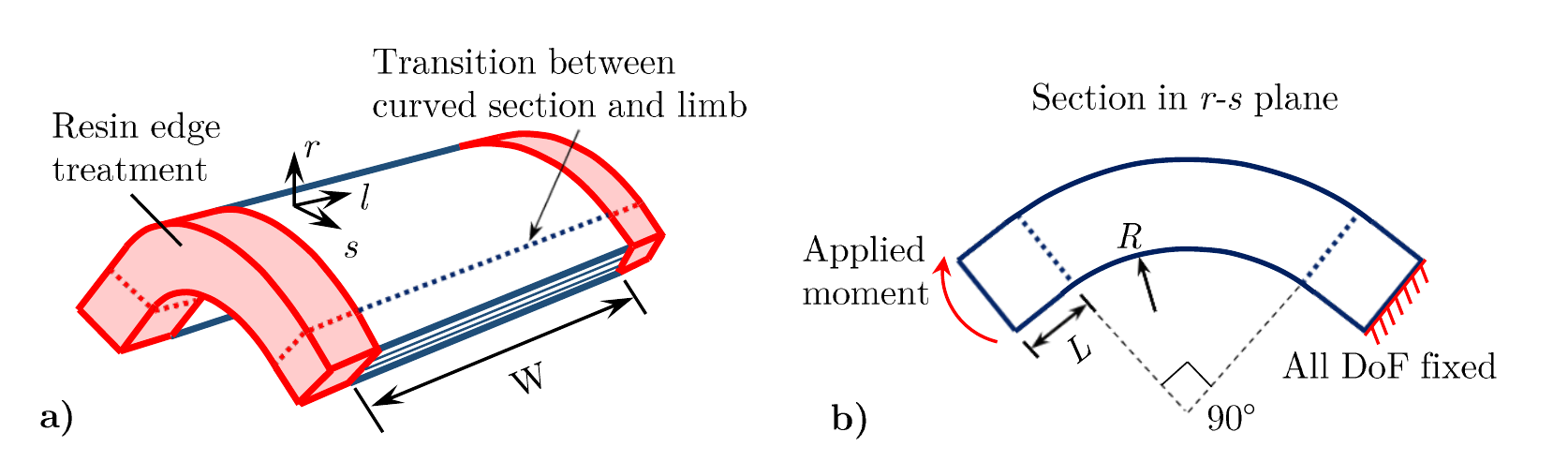}
 \caption{Sketch of the edge-treated corner-unfolding specimen (left) and a 
 cross section of the configuration (right).}
 \label{fig-mpcsetup}
\end{figure}

Figure \ref{fig-mpcsetup} shows the setup of our problem. For this small test we use a width
$W = $ 15mm and $L = $ 3mm long limbs. The inner radius $R$ of the curved section is 6.6mm.
The overall thickness $T$ is 2.98mm, where the ply layers are 0.23 mm thick and the resin
interfaces are 0.02 mm thick.
We then apply a resin treatment of 2mm to the free edges of the
laminate. This type of edge treatment has been shown in \cite{fletcher2016resin} to be
advantageous both for reducing conservatism in the design of aircraft
structures, as well as to make FE analyses more reliable.
The mechanical properties for the ply material and for the resin
interfaces of the curved laminate are given in Table \ref{table-1}.  
The fibre angles are given by the following $12$-ply stacking sequence
$$\left[\pm45/90/0/\mp45/\mp45/0/90/\pm45\right].$$

\begin{table}[t]
\centering
\begin{floatrow}
\begin{tabular}{ |c|c|c|c| } 
 \hline
 \multicolumn{2}{|c|}{\bf Orthotropic fibrous layer} &  \multicolumn{2}{|c|}{\bf Isotropic interface layer}\\
 \hline 
 $E_{11}$ &  162 GPa & $E$ & 10 GPa \\ 
 \hline
  $E_{22},\,E_{33}$ &  10 GPa & $\nu$ & 0.35 \\  
\hline
  $G_{12},\,G_{13}$ &  5.2 GPa & \multicolumn{2}{|c|}{}\\
 \hline
  $G_{23}$ &  3.5 GPa & \multicolumn{2}{|c|}{\bf Resin edge material}\\ 
 \hline
  $\nu_{12},\,\nu_{13}$ &  0.35 & $E$ & 8.5 GPa\\ 
 \hline
  $\nu_{23}$ &  0.5 & $\nu$ & 0.35 \\ 
 \hline
\end{tabular}
\caption{Mechanical properties for the curved $12$-ply laminate.}
\label{table-1} 
\end{floatrow}
\end{table}

In ABAQUS, the corner unfolding is modelled using a beam-type multi-point constraint (MPC) 
at the end of one limb. In DUNE,
  a similar effect is 
  achieved by adding a thin layer of very stiff material at the end of
  the same limb. In both cases, a $96.8$ Nmm/mm moment is applied to
  that limb, while fixing all degrees of freedom on the other limb, 
as shown in Figure \ref{fig-mpcsetup}.
We denote the stresses using the following directional notation:
arc length $s$ (around the curve); radius $r$ through-thickness and width 
$l$ across the laminate. These are shown by the coordinate system in Figure 
\ref{fig-mpcsetup}a.
  
DUNE and ABAQUS both use a model with second-order serendipity elements
(element \verb=C3D20R= in ABAQUS, \cite{serendipity}). There is a small difference in that
ABAQUS uses reduced integration while DUNE uses full integration.
Also, the stresses are not recovered in an identical way from the
displacements in the two codes. 
We use a tensor-product grid with $56$ elements in the $l$ and $s$ directions
and $6$ elements in each of the 12 fibrous and the 11 interface layers in
the $r$ direction. In total, this mesh contains $432,768$ elements. 
To ensure the effects at the free edge and at the material
discontinuities are sufficiently resolved, the FE mesh is graded
towards both edges in the $l$ direction, as well as towards the
interfaces in each of the resin and fibrous layers in the $r$
direction. The bias ratio 
gives the ratio between width of the smallest and largest element of the mesh.
We choose a bias ratio of $400$ between the
centre and the edges of the structure in the $l$ direction and a bias
ratio of $10$ between the centre and the interfaces in each layer in
the $r$ direction.

\begin{figure}[p]
\includegraphics[scale=0.93]{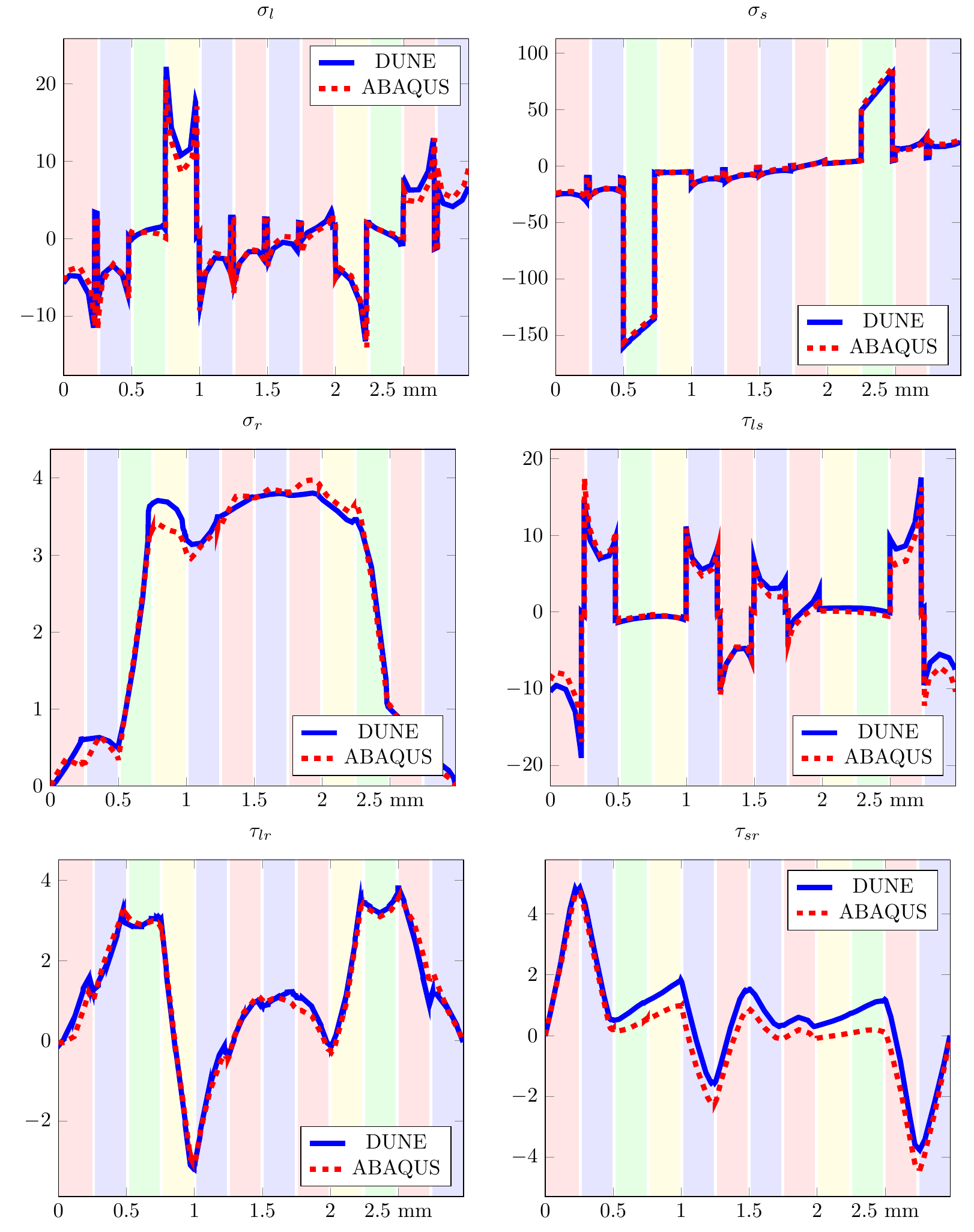}
   \caption{Stresses in MPa as a function of radius $r$ from outer radius
   at the apex of the curve,
   at $2.156$mm from the edge of the resin-edge-treated laminate
    (DUNE, solid blue; ABAQUS, dotted red).
    The background colours indicate the stacking sequence: $+45$ =
    red, $-45$ = blue,
    $90$ = green, $0$ = yellow.}
    \label{fig-stresses}
\end{figure}
 
In Figure \ref{fig-stresses},  all six stresses computed with ABAQUS and DUNE
are compared at the apex of the curve (i.e., at an angle of $45^\circ$) near the resin edge, $2.156$mm
away from the edge, which is in the laminate $0.156$mm ($6$ elements)
away from the resin laminate edge. Since we use identical meshes and identical
elements, the numbers of degrees of freedom are also the same. However, there are some differences
in the models, such as the quadrature order, use of single precision in ABAQUS and the stress
reconstruction, which could account for the different results. Nevertheless, the results 
show reasonably good agreement for all six different stresses.

\subsection{Efficiency of \textit{dune-composites}: Robust iterative solvers}\label{sec-solvers}
The computational gains of \textit{dune-composites} demonstrated in this
section are a direct consequence of the use of fast and robust
iterative solvers for the resulting large systems of linear equations
 \begin{equation}\label{eq-1}
 K\underline{u}=\underline{f},
 \end{equation}
where $K$ is the global stiffness matrix, $\underline f$ is the load
vector arising from the applied boundary conditions and $\underline u$
is the solution vector containing
the displacement values at each node within the finite element
mesh. There are also some significant gains in the setup time due 
to a more efficient, problem-adapted data management.

Solvers for systems of linear algebraic equations can broadly be
classified into direct and iterative ones. Direct solvers, based on
matrix factorisation, are more universally applicable and more robust
to ill-conditioning \cite{DirectSparse}. Thus, they are typically the default in
commercial FE packages. However, the cost of the factorisation and the
memory requirements can quickly become infeasibly high for large
problems. Moreover, it is difficult
to achieve good parallel efficiencies on large multicore computer
architectures. Since they do not require any factorisations, iterative 
methods for \eqref{eq-1} have the potential to scale
optimally, both with respect to problem size and with respect to the
number of processors in a parallel implementation, but that crucially
depends on the 'conditioning' of the problem.

In size, we refer to the total number of degrees of freedom $N$ in the
underlying finite element (FE) solution $\underline u$. Below we will see examples,
for which $N$ can be very large, up to 60 million. It is important to
note that $K$ is {\em sparse}, that is, the number of nonzero entries in $K$
is proportional to $N$, which is significantly less than the possible
$N^2$ entries in a full matrix. Iterative methods only require the
storage of the original, sparse stiffness matrix, while the storage of
the factors in direct methods is closer to the $N^2$ entries in a full
matrix. The computational cost of the best direct methods for FE systems
arising from three-dimensional structural calculations, such as those
used in ABAQUS, still grows at least with order $N^{1.5}$, as we will
see below. Iterative solvers that only require multiplication with the
sparse matrix $K$ have the potential to scale linearly with the problem size $N$.
 
The conditioning is
a property of the global stiffness matrix $K$. We say that $K$ is
ill-conditioned, if the nodal displacements $\underline{u}$ are
highly sensitive to rounding errors and to small changes in
$\underline f$. Generally, simulations of composite structures are
ill-conditioned because of the strong heterogeneity
and the often complex, non-grid aligned anisotropy of the material. 
The condition number of $K$, which is defined as the
ratio of the largest and the smallest eigenvalue of $K$, grows roughly
linearly with the size of the largest jump in the entries of the stiffness
tensor from one finite element to an adjacent one. 
But the conditioning also worsens systematically, as the problem size
$N$ increases. This growth is typically of order $N^{2/3}$ in 3D
problems.  This growth affects the accuracy of the solution in direct
solvers, but it has no effect on the computational
cost.

By contrast, for iterative methods, such as Gauss-Seidel, conjugate gradients
or GMRES methods, the conditioning of $K$ has a direct impact on the
number of iterations and thus on the cost. The basic idea of iterative
methods is to approximate the solution iteratively, starting from some
initial guess, by
sequentially reducing the residual $$r=K\underline{u}-\underline{f}$$
at each iteration \cite{IterativeSparse}. A 'good' or fast iterative solver reduces this
error quickly, in a few iterations. For well-conditioned problems,
standard black-box iterative solvers (as offered by ABAQUS) work
well. But the number of iterations typically grows with the
square-root of the condition number of $K$. For the poorly
conditioned problems of interest in this paper, 
standard iterative methods will either converge very slowly or
completely fail to do so. 

A standard approach to improve the conditioning of such problems is to
apply a {\em preconditioner} $P$ to (\ref{eq-1}), typically a cheap approximation of
$K^{-1}$, and then to iteratively solve
\begin{equation}\label{eq-pre}
P K\underline{u}=P \underline{f}. 
\end{equation}
The better the approximation $P \approx K^{-1}$, the closer the system
matrix in  \eqref{eq-pre} is to the identity and thus to a
well-conditioned problem with condition number close to $1$. However,
as outlined above, computing the inverse of $K$ (e.g. via
factorisation) is too expensive. Ideally the cost of applying $P$ and
the memory requirements should again only grow linearly with~$N$. 
A good choice of $P$ is problem
dependent, and must be carefully designed and tuned, to obtain 
fast solvers which scale well for large problems over multiple cores.
However, some well-understood common design choices for FE
discretisations of elliptic partial differential equations exist. To
achieve independence (or at most logarithmic dependence) of the number of
iterations on the problem size $N$, it is paramount to apply
multilevel preconditioners, either in the form of multigrid
\cite{BrHeMc, Vass} or multilevel domain decomposition methods 
\cite{domaindecomp}. 
In addition, the memory demands of such iterative
solvers are small (proportional to $N$) and the work
can be easily distributed over multiple cores. 

Care has to be taken though to achieve
robustness with respect to heterogeneities and anisotropies in the 
material properties. Since black-box multilevel approaches do not scale 
indefinitely and are not robust, more tailor-made approaches are needed. The
flexibility to prescribe or implement such tailor-made
preconditioners is not available in commercial software, and therefore the
 development of codes like \textit{dune-composites} is pivotal in solving
 large-scale problems in composite structures.

 \begin{figure}[t]
 \centering
 \includegraphics[scale = 0.8]{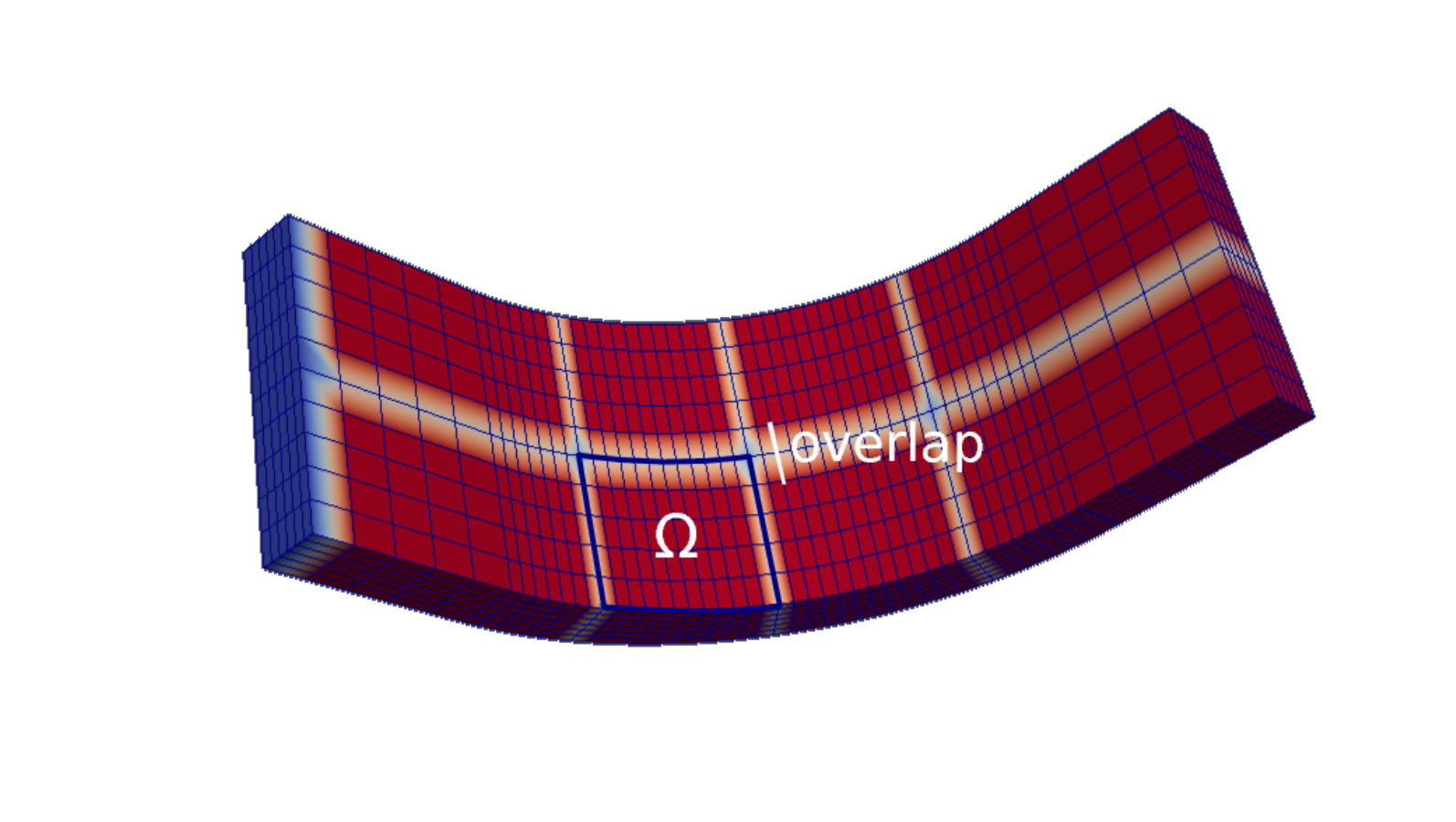}\\
 \caption{Decomposition of a mesh into several subdomains
   ($\Omega$), showing an overlap of one element in each subdomain.}
   \label{fig-domaindecomp}
   \end{figure}
 
In \textit{dune-composites} we use the conjugate gradient (CG) iterative
method for \eqref{eq-pre} and we test two types of preconditioner: the
generic Algebraic Multi Grid (AMG) preconditioner \cite{dune3}, as
well as a new two-level overlapping domain decomposition
preconditioner (GenEO) with a
multiscale coarse space based on local, generalised eigensolves in
each of the subdomains~\cite{geneo2}. 
For their parallelisation, both methods partition the physical 
domain and the FE mesh into subdomains of roughly equal size and distribute
them across multiple cores (see Figure \ref{fig-domaindecomp}). 
The key ingredient in both preconditioners is a coarse approximation
of $\underline{u}$ with significantly fewer degrees of freedom that is
cheaper to obtain then $\underline{u}$. This coarse approximation 
is then interpolated back to the original FE mesh. 

In Algebraic Multigrid, the problem is coarsened repeatedly by
aggregating degrees of freedom in the stiffness matrix $K$, creating a hierarchy of
coarse approximations of decreasing size. The aggregation process is
purely algebraic and based on the fact that the displacement
$u_i$ and $u_j$ at two neighbouring nodes $x_i$ and $x_j$ in the 
FE mesh will be similar if the nodes are ``strongly connected'', i.e. 
$u_i \approx u_j$ if $|K_{i,j}|$ is
large (relative to $|K_{ii}|$ and $|K_{jj}|$). The degrees of freedom
at two (or more) nodes that are strongly connected are then aggregated 
to a single degree freedom to obtain a cheaper approximation. This
process can be applied recursively until a sufficiently coarse
approximation has been obtained that can be efficiently solved using a
direct method. On the intermediate levels and on the original mesh we
only apply a cheap relaxation method such as Gauss-Seidel. The success of the
AMG preconditioner largely depends on how well this aggregation
process works (for details see \cite{dune3,Vass}). 

\begin{figure}[t]
 \centering
 \includegraphics[scale=0.14]{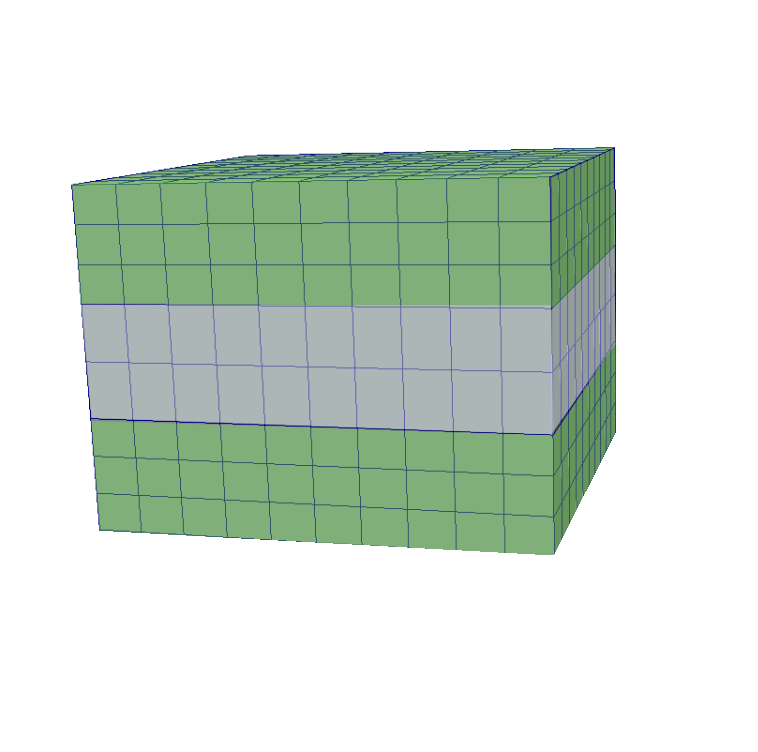}
 \includegraphics[scale=0.14]{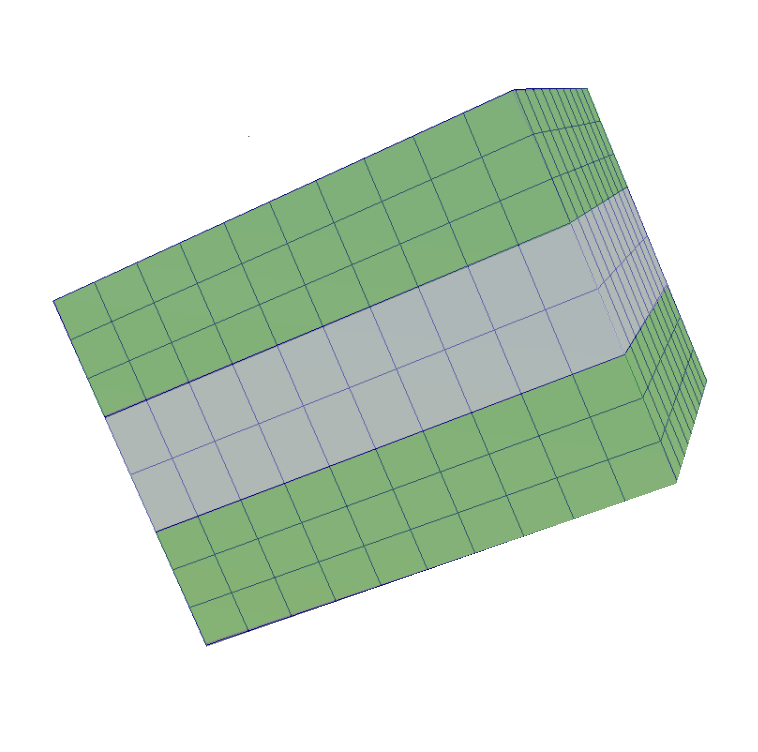}
 \includegraphics[scale=0.14]{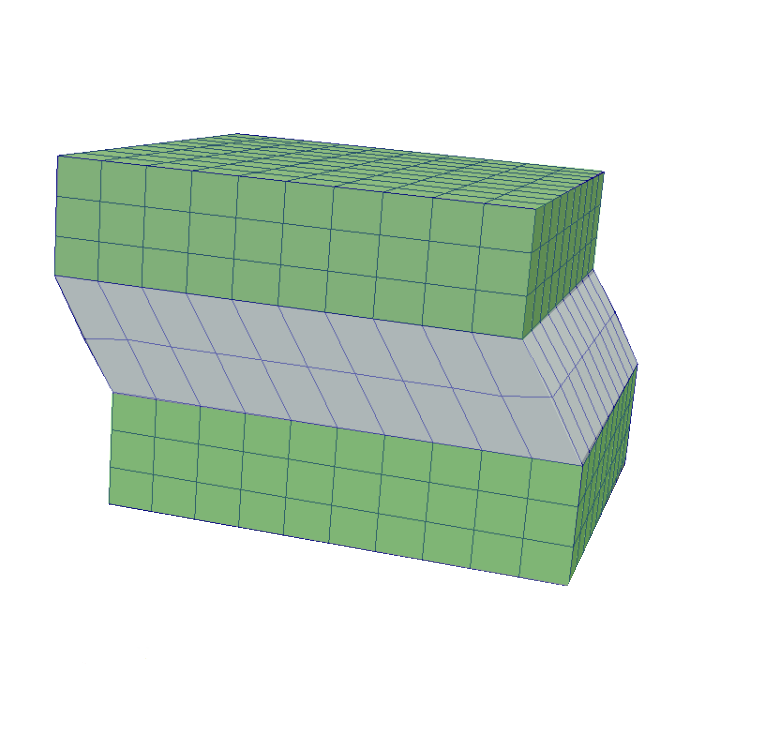}
  \includegraphics[scale=0.14]{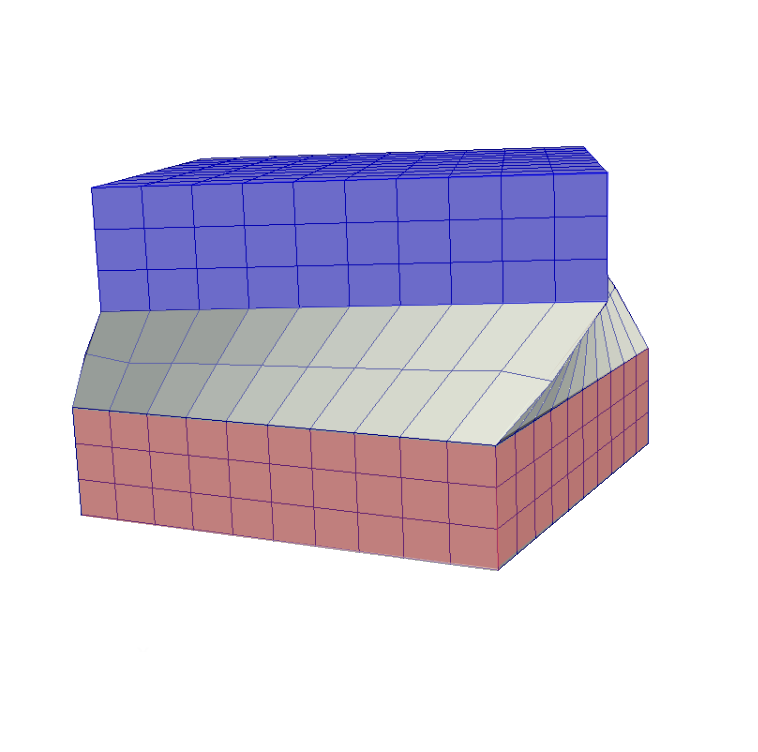} 
 \caption{Schematic of low energy eigenmodes in a reference block consisting of a layer of resin (white)
 between two layers of plies (green: $90^\circ$, blue: $-45^\circ$, red:
 $45^\circ$); from left to right, a reference block, 
 a zero energy mode (rotation),  a low energy mode
 between two $90^\circ$ plies and a low-energy twisting mode caused by $45^\circ$ and
 $-45^\circ$ plies.}
 \label{fig-eig}
\end{figure}

The GenEO preconditioner is of overlapping Schwarz-type
\cite{domaindecomp} and uses only two levels. The domain is
partitioned into overlapping subdomains of roughly equal size (see
Figure \ref{fig-domaindecomp}). In each of these subdomains, the 
original FE problem subject to zero-displacement boundary conditions 
on the artificial subdomain boundaries is solved either exactly using
a direct method or approximately using, for example the AMG
preconditioner. To capture the global behaviour of the solution, a coarse
problem is introduced. For homogeneous structures, standard finite 
elements on a coarser mesh are sufficient \cite{domaindecomp}. For composite
structures, especially in the presence of defects, this approach is
too inaccurate. Instead, the GenEO preconditioner computes a few of
the smallest energy (eigen) modes of the structure on each of the
subdomains. This includes the rigid body modes, but also some more exotic
modes that arise due to the particular properties of the composite. 
Some of these low energy modes can be seen in Figure
\ref{fig-eig}. Finally, these local modes are combined to 
generate a coarse global approximation using a so-called {\em
  partition of unity} approach \cite{domaindecomp}. The coarse problem is again solved
using a direct method or an iterative method. For details see the original paper
\cite{geneo2} and the related works \cite{EfGa, BabLip}. 
The GenEO coarse space has been shown to lead to a fully
robust two-level Schwarz preconditioner which scales well over
multiple cores \cite{geneo2,geneoparallel}. For isotropic problems
this has been proved rigorously in \cite{geneo2}. The robustness is
due to its good approximation
properties for problems with highly heterogeneous material
parameters. This is in fact of independent interest
\cite{BabLip,Dodwell17}.  Details of our specific implementation of GenEO in DUNE are given in
\cite{Seel}.  

\subsection{Performance Tests and Parallel Efficiency}\label{sec-time}

\begin{figure}
\centering
\includegraphics[scale=0.9]{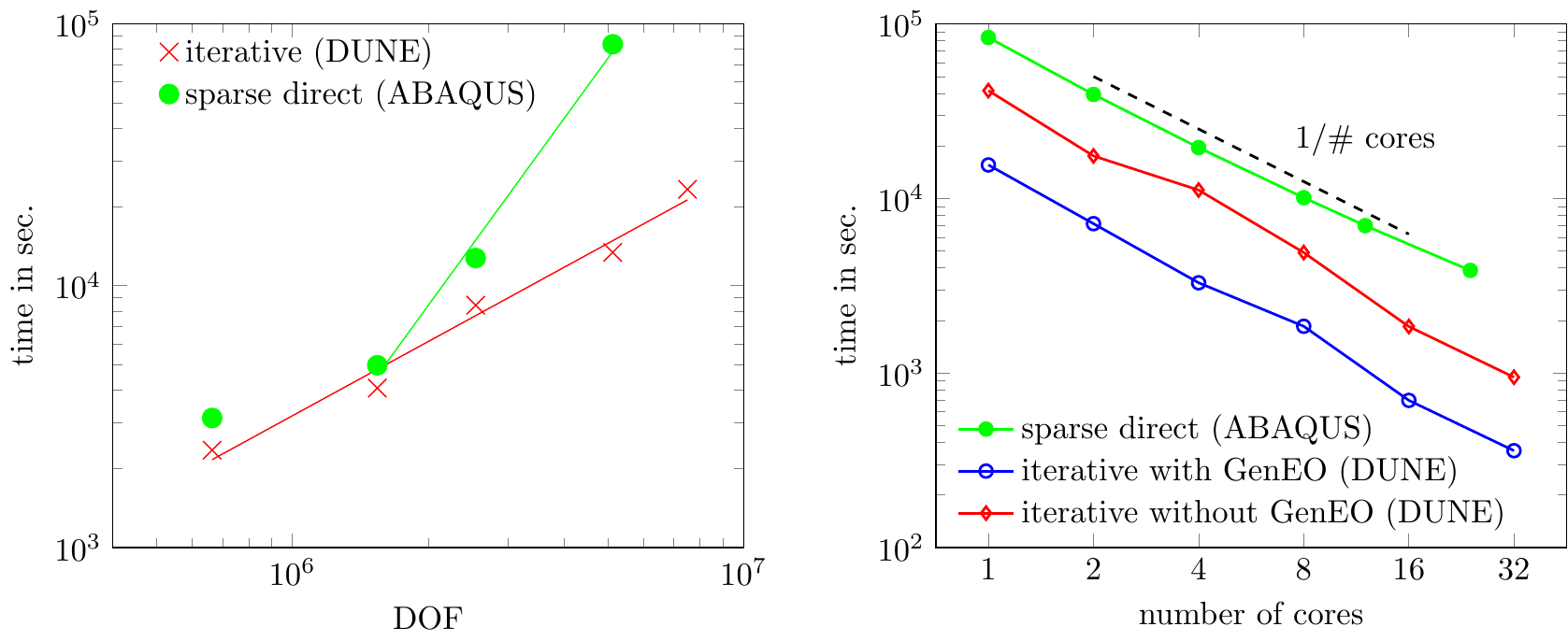}
 \caption{Left: Timing comparison in terms of degrees of freedom between the
   sparse direct solver in ABAQUS and the iterative CG solver (with AMG
   preconditioner) in DUNE. Right: Timing comparison in terms of numbers of cores between the
   sparse direct solver in ABAQUS and the iterative CG solver (with
   one- and two-level domain decomposition preconditioner) in DUNE.}
 \label{fig-sequential}
\end{figure}

First, we compare the performance of various preconditioners for the iterative solver in
\textit{dune-composites} when simulating the problem specified in Section
\ref{sec-verify}, as well as benchmarking them against the direct solver
in ABAQUS. The iterative solvers in ABAQUS only converge for
homogeneous problems or for very small numbers of degrees of freedom.
The measures for comparing the solvers are their overall computational
cost and their parallel efficiency. All the experiments in Figure
\ref{fig-sequential} were carried out on the computer {\tt cts04}
which consists of four 8-core Intel Xeon E5-4627v2 Ivybridge
processors, each running at 1.2 GHz and giving a total of 32 available cores.

On the left in Figure \ref{fig-sequential}, we compare the iterative solver with
AMG preconditioner in DUNE (with tolerance $10^{-4}$) and the sparse direct solver in
ABAQUS for increasing problem sizes. We see that the
iterative solver in DUNE scales linearly with
the problem size $N$ (the number of degrees of freedom) while the
cost of the direct solver in ABAQUS
grows significantly faster. Asymptotically, the growth is more than
quadratic in $N$. Moreover, the direct solver runs out of memory on one core for the 
biggest problem tested in DUNE ($N \approx 7.5$ million).
We specifically use the block version of
AMG in our tests in \textit{dune-composites} and
measure strength of connection 
between two blocks for the aggregation procedure in 
the Frobenius norm. As the smoother, we use two iterations of symmetric
successive overrelaxation (SSOR). The number of
iterations remains almost constant and moderate in this test. 
However, AMG has not originally been
designed for serendipity elements, and in its current form it does not
seem to be robust for this element,
especially in parallel.

The second comparison, on the right in Figure \ref{fig-sequential}, shows the
parallel scaling of various solvers. Here, we plot the strong scaling,
that is, the reduction in computational cost due to an increase in
the number of cores for a fixed problem size. We see that up to $32$
cores both the
parallel direct solver in
ABAQUS (green curve), as well as the iterative solver with domain
decomposition preconditioner in DUNE (red and blue curves) show almost 
optimal parallel scaling, that is, the CPU time is halved every time
the number of cores doubles. We tested two versions of the domain
decomposition preconditioner. The red line corresponds to a one-level
overlapping Schwarz preconditioner with an overlap of two layers of
finite elements. The subdomain problems are solved with UMFPACK \cite{umfpack}, a 
sparse direct solver, similar to the one in ABAQUS. For 1-8 cores, we use 8
subdomains and partition them evenly onto the available cores. For 16
and 32 cores, we use 16 and 32 subdomains, respectively (one per
core). The blue line
corresponds to the new GenEO preconditioner, a two-level overlapping 
Schwarz preconditioner where, in addition to the subdomain solves, we
also build and solve a problem-adapted coarse problem (see Section
\ref{sec-solvers}) . This coarse 
problem is again solved using UMFPACK. 

The addition of the domain decomposition preconditioner was essential 
for a good parallel efficiency. For some reasons which we cannot fully
explain, the parallel AMG preconditioner, which is available as a default within
DUNE did not scale well in parallel. The addition of the GenEO coarse space
dramatically reduces the number of iterations from about $400$ for the 
one-level overlapping Schwarz preconditioner (diamond points) to about $40$
for GenEO (hollow circle). Both those numbers remain roughly
constant across this range of cores. In absolute terms, the lower
number of iterations for GenEO translates into an almost three-fold 
reduction in computational time. Moreover, beyond 32 cores, only
the GenEO preconditioner continues to scale. The numbers of iterations
for the one-level Schwarz preconditioner start to increase with the
number of subdomains.

Finally, we note that even for the relatively small
test we have chosen here ($432,768$ elements) the iterative solver in DUNE
preconditioned with GenEO is around $5-10$ times faster than the direct 
solver in ABAQUS. Moreover, the restriction to shared-memory
parallelism and the memory requirements of the direct solver in
ABAQUS, combined with the $O(N^{2})$ growth in computational cost,
make it essentially impossible to carry out a corresponding scaling
test with ABAQUS beyond $32$ cores.

\begin{figure}[t]
 \centering
 \includegraphics[scale=1.1]{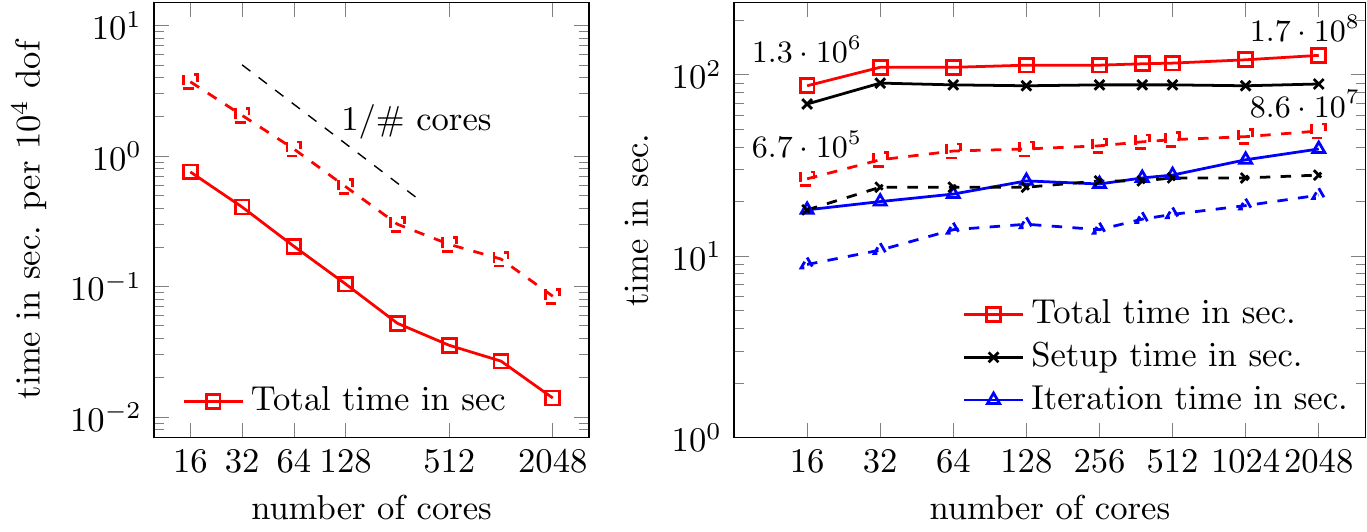}
 \caption{Scaling of {\em dune-composites} using the GenEO
   preconditioner on up to $2048$ cores. Solid lines are used for the problem with
 $8.4\cdot10^4$ degrees of freedom per core; dashed lines for 
$4.2\cdot10^4$ degrees of freedom per core. For the larger problem,
the problem size varies from about $1.3$ million degrees of freedom on
16 cores to about $173$ million degrees of freedom on 2048 cores.}
\label{fig-geneo}
\end{figure}

By contrast, in Figure \ref{fig-geneo}, we now look at the parallel
scaling of the iterative solver with GenEO preconditioner in
\textit{dune-composites} on thousands of cores solving problems with
up to $173$ million degrees of freedom. We use again UMFPACK for all 
subdomain and coarse solves in GenEO. 
For this second scaling test, we move to a larger scale
problem. 
The model run in this section has $39$ plies and a variable
width. For details on the setup of this
problem see Section \ref{sec-def}. We also move to a larger computer. 
All the experiments in Figure \ref{fig-geneo} were carried out on the 
University of Bath HPC cluster {\tt Balena}. This consists of 
192 nodes each with two 8-core Intel Xeon E5-2650v2 Ivybridge
processors, each running at 2.6 GHz and giving a total of 3072 available cores. 

We carry out a weak scaling
experiment, that is, we increase the problem size proportionally to
the number of cores used. For an iterative solver that scales optimally
both with respect to problem size and with the number of cores,
the computational time should remain constant in this experiment. 
To scale the problem size as the number of cores $N_\text{cores}$
grows, we increase the width of the laminate.  We present results for
two different setups. For the first, larger problem we have used $56$
elements along the radius, $\frac58 N_\text{cores}$ elements across the
width and $4$ elements through thickness in each of the resin and ply
layers, respectively. This setup ensures that the amount of work
handled by each core remains constant as $N_\text{cores}$ varies.
For the second, smaller problem we have halved the number
of elements through the thickness, using only 2 elements in 
each resin and ply layer. In Figure \ref{fig-geneo}, we see that after
a slight initial growth the 
scaling of the iterative solver in DUNE with GenEO preconditioner
is indeed almost optimal to at least $2048$ cores, allowing us to 
increase the size of the tests at a nearly constant run time and thus, to solve a
problem with $173$ million degrees of freedom in just over 2 minutes.

 \section{Defect analysis and the influence of boundaries}\label{sec-def}
  
\begin{figure}[t]
 \begin{floatrow}
   \centering
  \includegraphics[scale=0.35]{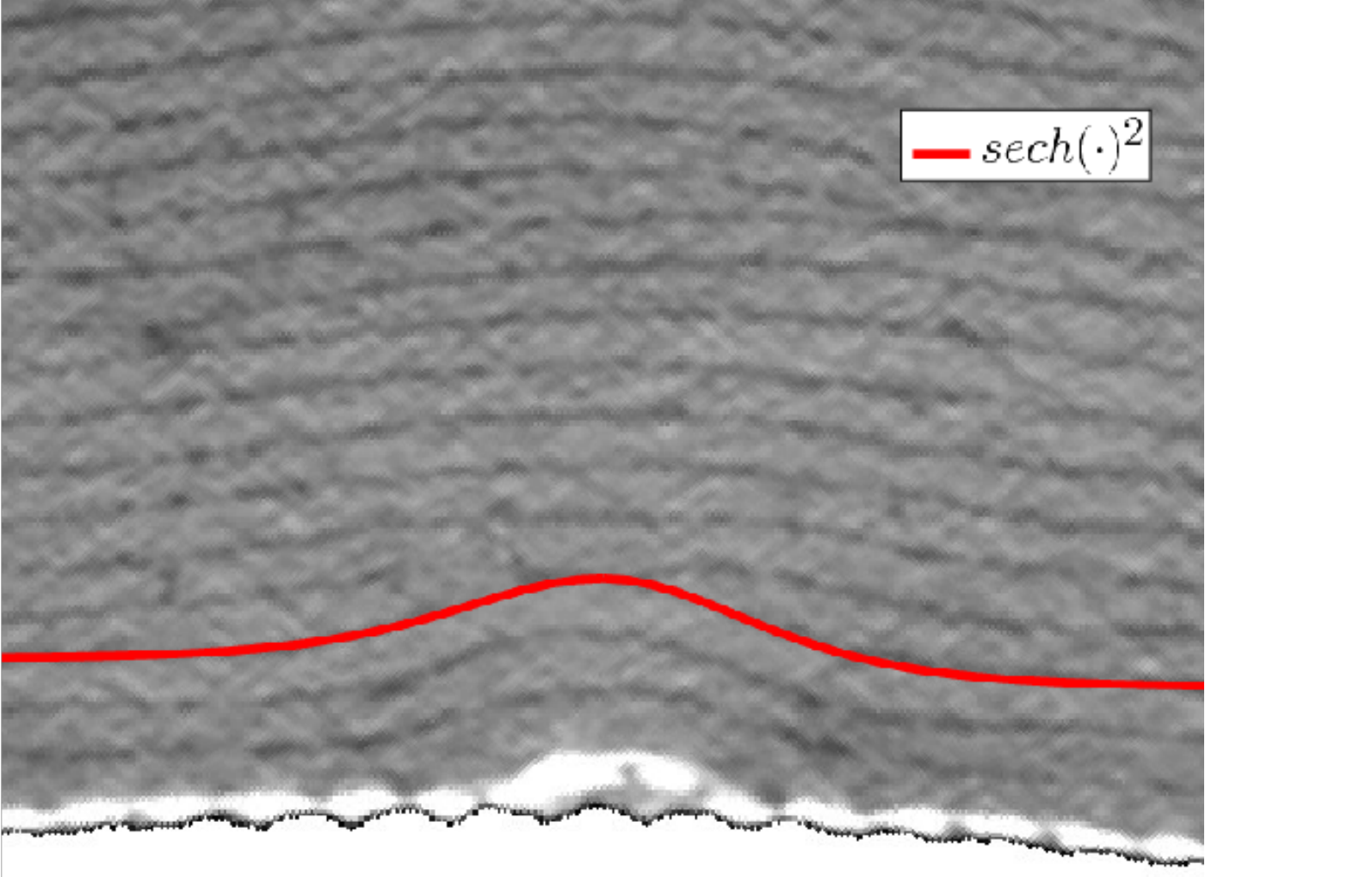}
  \includegraphics[scale=0.17]{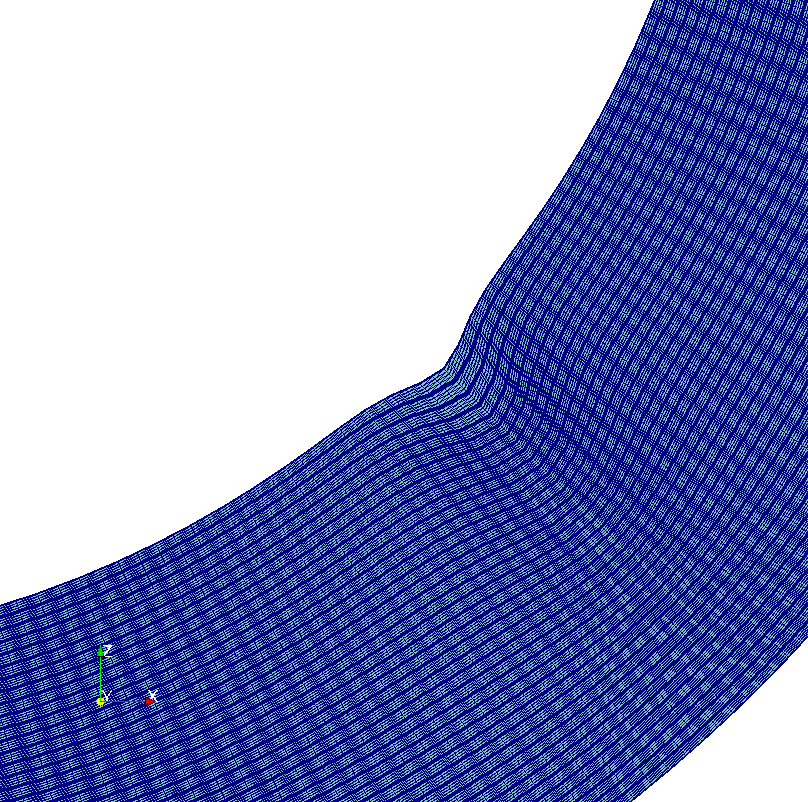}
  \includegraphics[scale=0.17]{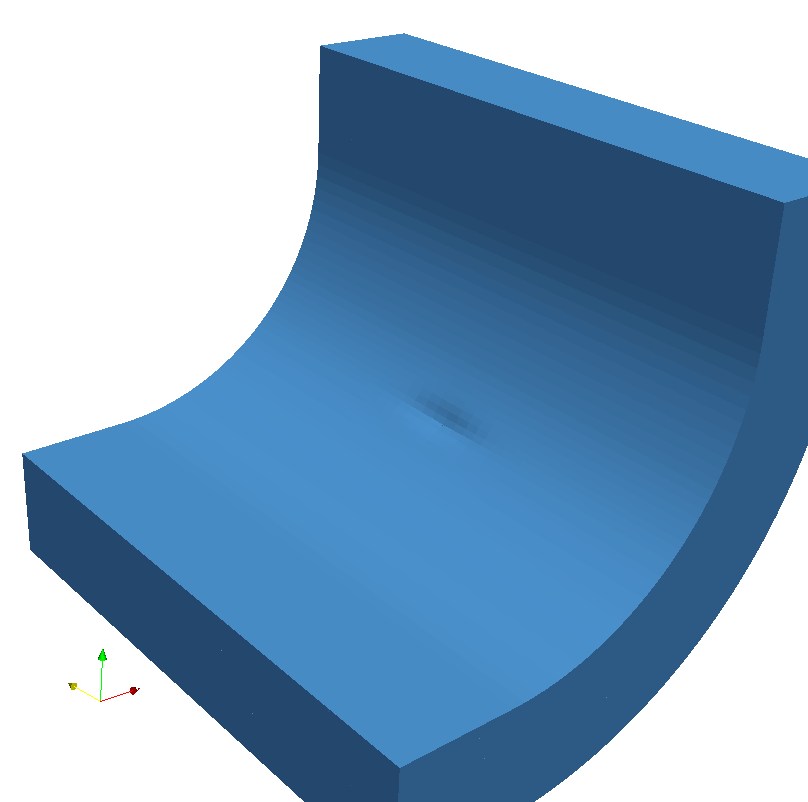}
 \end{floatrow}
   \caption{Image of a defect. Left: A scan of a wrinkle with the $\text{sech}^2\,(\cdot)$ superimposed.
   Center: A radial cut   through the center of the laminate. 
   Right: Decay of the wrinkle across the width.}
  \label{fig-def}
 \end{figure}

 In manufacturing, small localised defects in the form of misaligned fibrous layers can occur,
 these defects can have a large effect on the strength of the materials. 
 In this section we investigate  the effect that varying the maximum slope and amplitude of 
 a localised wrinkle with only one oscillation has on the strength of the curved laminate.
 We show that we need a very high-fidelity mesh to be able to compute localised stresses and this
 problem provides an appropriate application for \textit{dune-composites}.
In order to determine the shape and size of the wrinkle we used an  X-ray CT scan of a typical curved
laminate with a wrinkle defect, shown in  Figure \ref{fig-def} and described in more detail
below.

 For the defect analysis we we examine the influence of width-wise boundary conditions on the results.
 The mechanical properties for the ply material and the resin interfaces
 are the same as those for the previous smaller test, see Table \ref{table-1}.
 The geometry of this $39$ ply problem is chosen to give a similar ratio of
 width to radius as the previous test. Here we use a width
$W$ of 52mm and 10mm long limbs $L$. The inner radius $R$ of the curved section is 22mm.
The overall thickness $T$ is 9.93 mm, where the ply layers are 0.24 mm thick and the resin
interfaces are 0.015 mm thick.
 The fibre angles of plies are given by the following stacking sequence 
 $$\left[[\mp45/90/0]_2/[\mp45]_2/90/\mp45/90/0/\mp45/0/\pm45/0/90/\pm45/90/
 [\pm45]_2/[0/90/\pm45]_2\right].$$

   \begin{figure}
 \centering
 \includegraphics{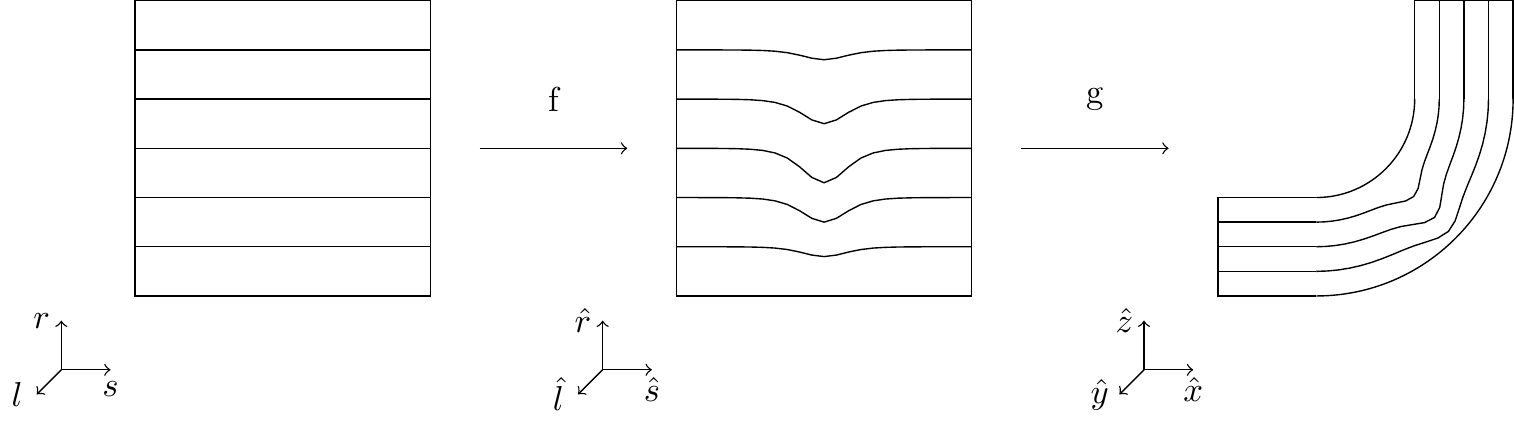}
\caption{Transformation of the flat geometry into a flat geometry with a defect and finally into
the curved geometry with a defect.}
\label{fig-transformation}
 \end{figure}

The boundary conditions and loading in a narrow specimen are very different
from those of the full structure. For this reason we initially model using 
periodic boundary conditions
at the free edge. This gives an approximation of the effect of the wrinkle in a very wide part.
 For the later tests we remove the periodic boundary conditions and add a $3$mm wide
 layer of resin to the free edges. This reduces the strength of the singularity
 at this edge and for larger wrinkles, ensures that failure occurs
 at the center of the laminate, near the defect.

Figure \ref{fig-transformation} shows the process by which the wrinkle is created.
We start with a mesh on a flat plate, then we add the wrinkle to the flat plate and
finally transform to the curved geometry. Let $(s,l,r)$ be coordinates in the unperturbed, 
flat geometry. Note, that in the final curved geometry, $s$ corresponds to an arc length, 
$r$ is the radius through the part
and $l$ is the width across the part.

In the flat geometry the wrinkle is given by a perturbation in $r$, which depends
on all three coordinates. The perturbed
coordinates $(\hat s, \hat r, \hat l)$ 
are given by
\begin{equation}\label{eq-wrinkle}
\hat s = s, \quad \hat l = l \quad \text{and} \quad \hat r = r+f(s,l,r).
\end{equation}
If we denote by $(s_\text{def},l_\text{def}, r_\text{def})$
the location where the defect is largest, the perturbation is given by
\begin{equation}\label{eq-f}
 f(s,l,r) = d \,\text{sech}^2\left( \pi \frac{s-s_\text{def}}{\pi b_1}\right)
                        \,\text{sech}^2\left( \pi \frac{r-r_\text{def}}{R b_2}\right)
                       \,\text{sech}^2\left( \pi \frac{l-l_\text{def}}{W b_3}\right),
\end{equation}
where $d$ gives the amplitude of the defect and $b_i$, $i=1,2,3,$ are parameters 
giving the ``extent'' of the defect in the three directions.

For these tests we have chosen $b_2$ so that the wrinkle decays more quickly towards the inner radius,
more precisely
$$b_2 = \begin{cases}\frac{1}{2}  & \text{ if } r-r_\text{def}<0\\ 
\frac{1}{4} & \text{ if } r-r_\text{def}>0\end{cases}.$$ Further,
$b_1=\frac{1}{5}$ and $b_3=\frac{1}{2}$. As shown in Figure \ref{fig-conana}, the sech$^2(\cdot)$ 
function never reaches zero, however at $\pm b_i$ it is already negligibly small ($<10^{-2}$). The coordinate $s$ ranges from $0$ to $\pi$, $r$ ranges from $R$ to $R+T$ and $l$ ranges from 
$0$ to $W$, which implies that a value of $b_3=\frac{1}{2}$ leads to a wrinkle spanning half of the width of the part.

The wrinkled mesh on the cube is finally mapped to the curved section via the mapping
\begin{equation}
 \begin{split}
  \hat x &= g_1(\hat s, \hat l,\hat r) = L + \hat r\sin\left(\frac{\pi}{2} \hat s\right)\\
  \hat y &= g_2(\hat s,\hat l,\hat r)  = \hat l\\
  \hat z &= g_3(\hat s,\hat l,\hat r)  = \hat r \cos\left(\frac{\pi}{2} \hat s\right),
 \end{split}
\end{equation}
 where
$L$ is the length of the limb. 

 From the form of the wrinkle we can easily calculate the steepest slope by deriving 
 \begin{equation}
  \tan^{-1}\left(\frac{d}{ds} f(s,l,r)\right) = \tan^{-1}\left(-\frac{2}{b_1}\,
  \text{tanh}\left( \pi \frac{s-s_\text{def}}{\pi b_1}\right) f(s,l,r)\right)
 \end{equation}
 This allows us to classify wrinkles according to the steepest angle produced by the defect. 
In order to ensure that the edge effect does not affect the results we give results with
periodic boundary conditions.
In Figure \ref{fig-def} we give a sample defect from a  CT scan with an overlay showing
the fit to our proposed parameterisation of the wrinkle. We also show the resulting model with
the largest amplitude (and thus slope) of defect considered in the following numerical tests.

\subsection{Convergence Analysis}
 The FE model used for these tests consists of approximately $7.3\cdot10^6$ degrees of freedom 
 in total. There are $96$ elements along the radius, $40$ across the width and
 $4$ each in the resin and ply layers. To resolve the defect the model is refined towards 
 the location of the defect, with a bias ratio of $2$, i.e. the elements
 at the center of the model are half as large as those at the edges. For the edge
 treated case there is also refinement towards the resin edge with a bias ratio of
 $400$ to resolve the stress concentrations that occur near the edge.
 
\begin{figure}[tb]
 \begin{minipage}{22em}
    \includegraphics[scale=1]{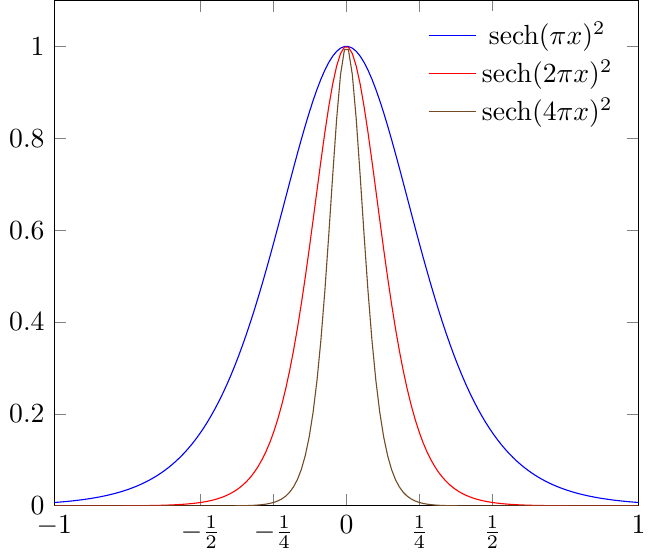}
 \end{minipage}
 \begin{minipage}{22em}
    \includegraphics[scale=.9]{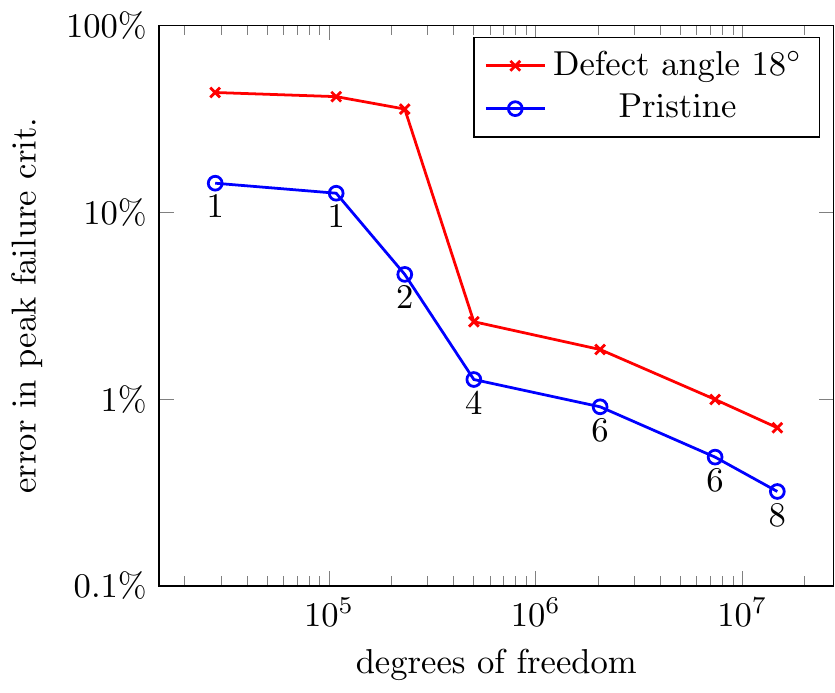}
 \end{minipage}
  \caption{Left: Decay of the sech$^2()$ function for $3$ different values of parameter $B$. 
  Right: Convergence plot for a pristine part and a part with a defect with maximum slope $18^\circ$, with
  labels giving the number of elements through-thickness per ply and interface layer.}
  \label{fig-conana}
 \end{figure}

 To ensure that the FE modelling is sufficiently accurate we include a convergence analysis.
 In Figure \ref{fig-conana} (right) we plot the relative error in the peak failure criterion given by
 equation (\ref{eq-Camanho}) below.
 As an approximation to the exact solution we used the next refinement level of the FE model
 ($136$ elements along the radius, $60$ across the width and $8$ each in the resin and ply layers).
 For the refinement we have chosen (second to last point) the relative error has been reduced
 to around $1\%$ for both the pristine model and for the model containing a defect.
 For the case with a wrinkle it is important to note that in the
 first three refinement levels have very similar peak failure criterion, even though the model
 is not yet converged. However,
 in these configurations the wrinkle is not yet sufficiently resolved and the failure criterion
 output is still closer to that of the pristine model. 
 
 The models used in this analysis increase the number of elements through thickness
 in each layer as follows: $1$, $1$, $2$, $4$, $6$, $6$, $8$.
The number of elements across the part and around the radius and limbs
are increased simultaneously. We can see that (especially for a part with
 a defect) the number of elements through thickness is very important. With fewer than 
 $4$ elements we do not get a good approximation of the maximum failure criterion.
 Further, when fewer elements are used the stresses drop to zero at the surface
 of the laminate.

 Around the curve enough elements need to be used to ensure that the wrinkle
 is resolved, at least $10$ elements should be used along the length of the wrinkle. Even
 with a mesh grading towards the wrinkle this still requires a large total number of elements
 across the full length of the part.

\subsection{Defect analysis} 
The stresses in the vicinity of a defect are complex, including
interlaminar, shear and direct stresses. Therefore a mixed mode failure criterion
is more suitable than a maximum stress criterion. The
strength of the laminates was assessed using a quadratic damage
onset criterion, defined by Camanho et al. \cite{Camanho} as
\begin{equation}\label{eq-Camanho}
 \sqrt{\left(\frac{\sigma^+_{r}}{s_{33}}\right)^2+\left(\frac{\tau_{rs}}{s_{13}}\right)^2
 +\left(\frac{\tau_{rl}}{s_{13}}\right)^2}=F
\end{equation}
with negative values of $\sigma_{r}$ treated as zero and the following allowables:
\begin{center}
$s_{33}=$ 61 MPa \qquad $s_{13}=$  97 MPa.
\end{center}

\begin{figure}[tb]
\centering
\includegraphics[scale=0.4]{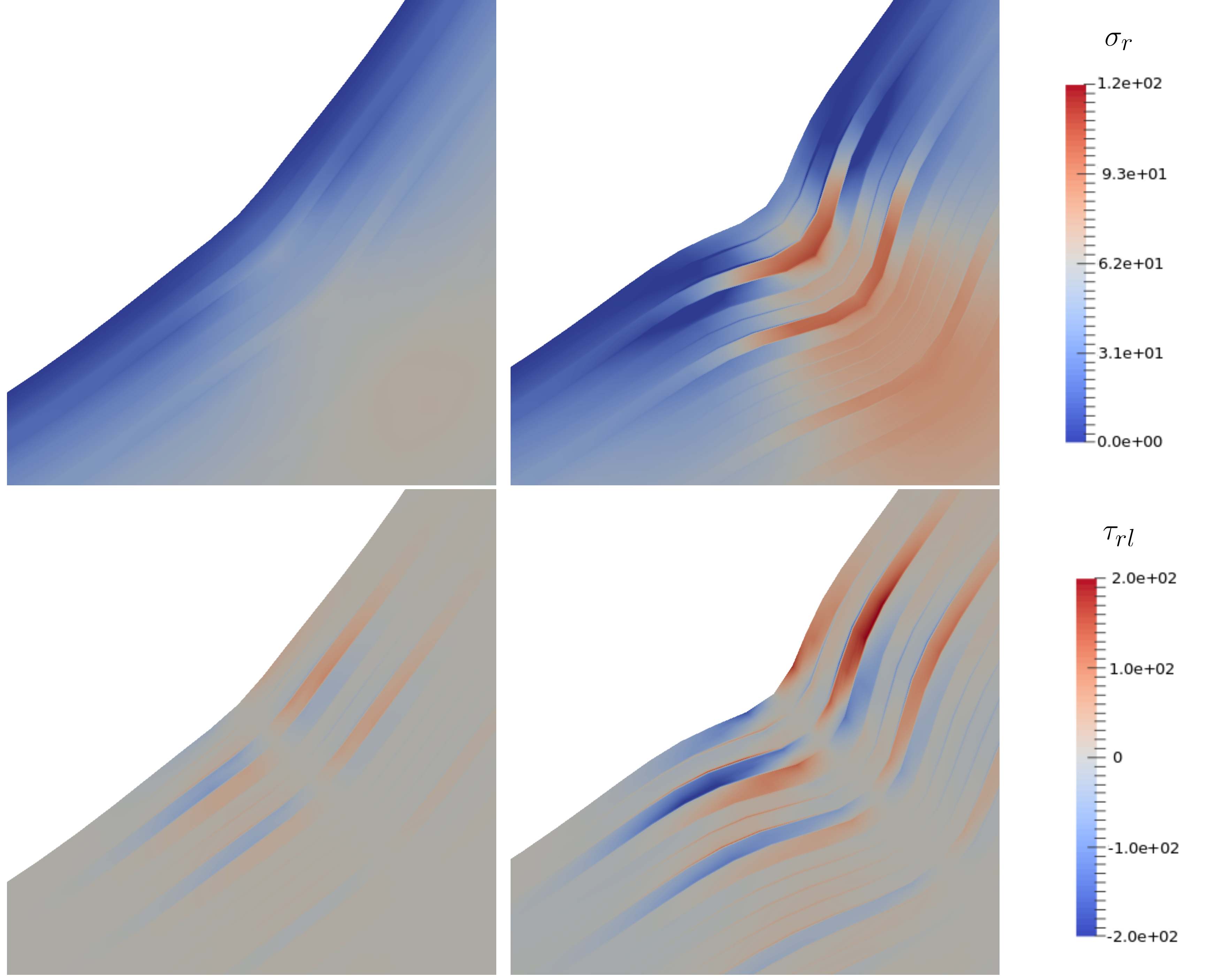} 
\caption{Stresses $\tau_{rl}$  and $\sigma_{l}$ for defect angle $4^\circ$  (left) , followed by 
stresses $\tau_{rl}$  and $\sigma_{l}$ for 
defect angle $18^\circ$ (right).}
\label{fig-stress}
\end{figure}

Curved laminates subjected to corner unfolding are observed to fail by delamination 
of the plies \cite{fletcher2016resin}.  Since this indicates failure occurs at the interface between plies, we 
apply the failure criterion only in the resin-rich interface zones and failure initiates when $F=1$.
Generally, the failure 
criterion is applicable in the local coordinates of the material, however since the interface
zone is isotropic, a transformation is not required.

\begin{figure}[H]
\centering
\includegraphics{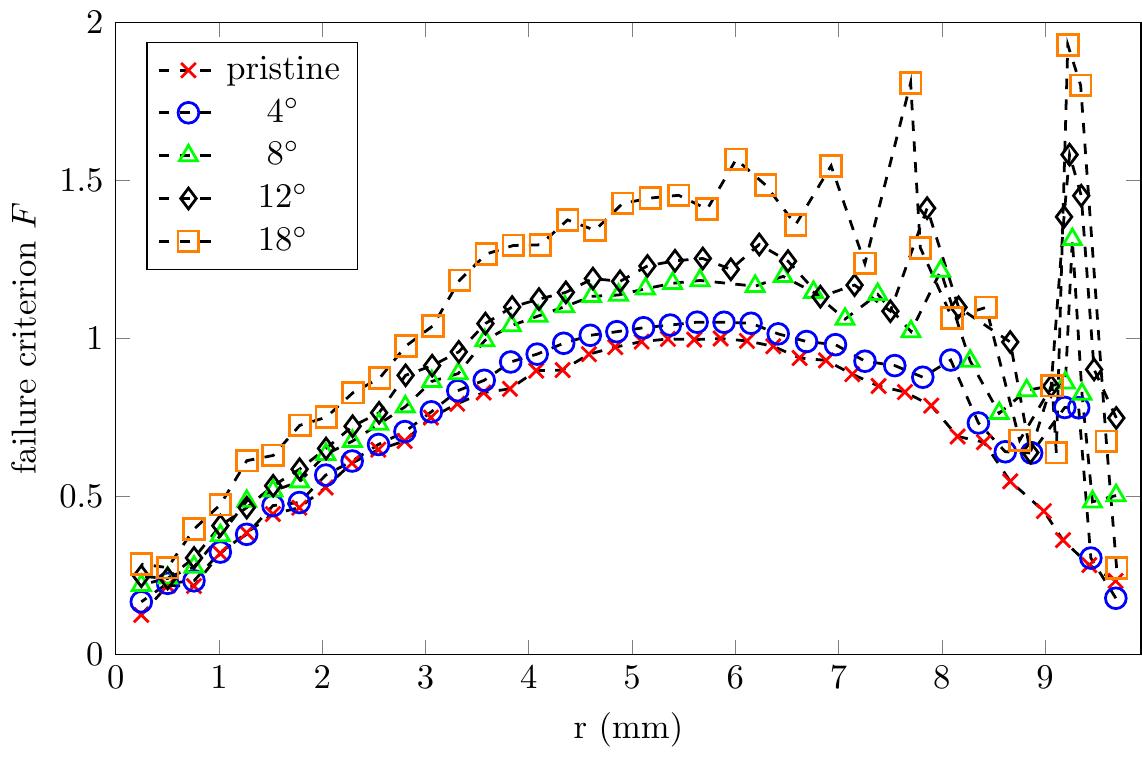}
  \caption{Through-thickness quadratic damage onset criterion $F$ at the center of each resin interface
  (with linear interpolation  between resin interfaces). Values are calculated at mid-width along 
  the angle of steepest slope ($2.1^\circ$ from  the apex of the curve), and at a distance $r$ from the
  outer radius.}
  \label{fig-max}
 \end{figure}

Having established the wrinkle defect model, we now investigate the effect of varying the wrinkle
severity.  Figure \ref{fig-stress} shows the effect on the interlaminar shear and direct stresses
with a wrinkle of maximum angle (slope) of $4^\circ$ and $18^\circ$ for which an opening/bending
moment of 9.58 kN mm/mm is applied.  The interlaminar shear 
maximises in the region of greatest slope, while the interlaminar direct stress maximises in
the region of peak amplitude.  These stresses are highly localised and cause a sharp rise in
failure index through-thickness, as shown in Figure \ref{fig-max}.  The pristine part is predicted 
to fail near the mid-thickness, where direct interlaminar stress maximises.  With a wrinkle 
of slope $4^\circ$, there is a rise in failure index near the inner radius, however failure 
is still predicted to first occur near the mid-thickness (as per the pristine model).
As the slope of the wrinkle is increased, so failure index near the inner radius becomes 
higher than that near the mid-thickness.  This is also indicated in Table \ref{table-per},
showing failure at interface $16$ for pristine and $4^\circ$ wrinkle models,
and at interface $4$ for greater slopes.  For the most severe wrinkle of slope $18^\circ$, 
the failure index is almost double that of the pristine model, and
hence strength is halved.

A summary of all models is shown in Table \ref{table-per}.  As seen by the maximum recorded
displacement, the overall stiffness of the curved laminate is not significantly affected by
the inclusion of a wrinkle defect.  This is because it is very small in comparison to the size
of the whole laminate.  However, it does significantly reduce strength, based on initiation of failure.
When tested in isolation such specimens have been found to fail in a highly unstable manner, suggesting once initiation
is reached there is instantaneous propagation and catastrophic failure \cite{fletcher2016resin}. 
Note that the two pristine results in Table \ref{table-per} compare well with predictions
acquired by modelling the curved laminates using ABAQUS software \cite{fletcher2016resin}. Note also that the edge-treated pristine result compares well with experimental results \cite{fletcher2016resin}. Table \ref{table-per} shows that the greatest impact of the wrinkle is to introduce significant 
interlaminar shear stresses.  Interlaminar direct stress is also increased, by approximately
$50\%$ between the pristine and $18^\circ$ wrinkle models; however, 
the interlaminar shear stresses go from being negligible in the pristine model,
to becoming greater than the interlaminar direct stress in the $18^\circ$ model.

Note that these results do not take into account thermal pre-stresses, which occur as a result of the
high-temperature curing process and may be significant. However, the purpose of this paper is to illustrate
the importance of mesh refinement for rapidly varying through-thickness stresses.

\section{Conclusions}
In this paper, a new high-performance FE analysis tool for composite
structures
is presented and shown to accurately model the 
behaviour of complex laminate structures with manufacturing
defects in a fraction of the computational time required by commercial
software packages, such as ABAQUS.
\begin{table}[H]
\resizebox{\textwidth}{!}{%
\begin{tabular}[scale=0.5]{ |p{3em}|p{3em}|p{3.4em}|p{3.4em}|p{3.4em}|c|c|c|c|c|c|p{3.4em}|} 
\hline
\multicolumn{12}{|c|}{\textbf{FEA (normalised): Periodic BC }}\\
\hline
 Defect slope (deg.) & max. disp. (mm) & location & interface (from outer radius) & max. failure crit. $F$
 & \begin{tabular}[c]{@{}c@{}}\\[-.5em]$\sigma_{l}$\\ (MPa)\end{tabular} 
 & \begin{tabular}[c]{@{}c@{}}\\[-.5em]$\sigma_s$\\ (MPa)\end{tabular}
 & \begin{tabular}[c]{@{}c@{}}\\[-.5em]$\sigma_r$ \\ (MPa) \end{tabular}
 &\begin{tabular}[c]{@{}c@{}}\\[-.5em] $\tau_{sl}$ \\ (MPa) \end{tabular}
 & \begin{tabular}[c]{@{}c@{}}\\[-.5em]$\tau_{rs}$ \\ (MPa) \end{tabular}
 & \begin{tabular}[c]{@{}c@{}}\\[-.5em]$\tau_{rl}$ \\ (MPa) \end{tabular}
 & M(kN)\\[-.2em]
 \hline
 0  & 3.99 & Mid- width & 23 & 1.00 & 33.76 & 34.90  & 61.00 & 0.19     & 0.02  & $-0.23$ & \textbf{9.58}$^1$  \\
 4  & 3.98 & Mid- width & 23 & 1.07 & 37.26 & 40.72  & 65.59 & 0.17     & 0.95  & 0.47    & 8.90  \\
 8  & 3.98 & Mid- width & 35  & 1.30 & 51.67 & 113.2 & 65.69 & $-8.75$  & 6.50  & 71.05   & 7.34  \\
 12 & 3.98 & Mid- width & 35  & 1.65 & 45.41 & 117.4 & 79.94 & $-12.80$ & 9.97  & 97.59   & 5.78  \\
 18 & 3.98 & Mid- width & 35  & 1.94 & 37.31 & 131.6 & 92.99 & $-17.37$ & 14.91 & 116.4  & 4.92  \\
 \hline
 \multicolumn{12}{|c|}{\textbf{FEA (normalised): Resin edge treated }}\\
\hline
 0   & 5.16 & Edge       & 23 & 1.16 & 80.06    & 116.25     & 11.13 & 6.75      & 28.13  & -107.5 & \textbf{8.26}$^2$ \\
 18  & 5.16 & Mid- width & 35 & 1.99 & 37.68    & 144.15     & 96.67 & $-16.72$  & 116.80 & 0.57   & 4.79 \\
 \hline
 \multicolumn{12}{r}{}\\[-.5em]
\multicolumn{12}{r}{\footnotesize{$^1$ ABAQUS result was $9.51$ kN \cite{fletcher2016resin}}}\\
\multicolumn{12}{r}{\footnotesize{$^2$ ABAQUS result was $8.25$ kN, test average $8.65$ kN\cite{fletcher2016resin}} }
 \end{tabular} %
 }
 \caption{Effect of the defect angle on peak stresses, on the failure criterion $F$ in interfaces
 and on M the predicted  moment of initiation of failure. }
 \label{table-per}
\end{table}

The stress field around a wrinkle defect in composite materials
is highly localised, requiring a high fidelity mesh to adequately
capture it.  Such defects that arise in composite materials during
manufacture are by nature difficult to predict.  There is a case
for performing stochastic analysis: running thousands of
    simulations with defects in various locations and sizes,
    studying the influence on the integrity of a component.
The high number of simulations, combined with the requirement for
high fidelity in each simulation, results in a very demanding
problem computationally.  The use of a bespoke FE solver, such 
as the one developed here in DUNE, is essential for meeting these demands.

Current industry standard FE tools, e.g. ABAQUS, are not
able to deal with these problem sizes, largely due to the direct
linear equation solvers that are employed and the restricted 
parallel scalability. Iterative solvers within ABAQUS only converge in
very simple cases. For this reason, as part of our project, a new 
iterative solver had to be implemented within DUNE. The heart of this
solver is a robust preconditioner (GenEO) that leads to an almost
optimal scaling of the computational time both with respect to the
number of degrees of freedom and the number of compute cores. This has
been demonstrated for problems with up to $173$ million degrees of
freedom, achieving run times of just over $2$ minutes on $2048$ cores.

The chosen  example problem evaluates the through-thickness stresses 
caused by unfolding a curved laminate.
This serves to illustrate the importance of using high-fidelity meshes
and provides insight into the
influence of small manufacturing induced defects.

 \section{Acknowledgements}
 This work was supported by an EPSRC Maths for Manufacturing grant (EP/K031368/1).
 Richard Butler holds a Royal Academy of Engineering-GKN Aerospace Research Chair in
 Composites. This research made use of the Balena High Performance Computing
 Service at the University of Bath.

\bibliographystyle{abbrv}
\bibliography{ElastPaper}

\end{document}